\documentclass[conference]{IEEEtran}

\IEEEoverridecommandlockouts

\usepackage{textcomp}

\usepackage{hyperref}

\usepackage{amsmath,amsthm,amsfonts,amssymb}
\usepackage[ruled,vlined]{algorithm2e}
\usepackage{xcolor}
\usepackage{enumerate}
\usepackage{cite}
\usepackage{csquotes}

\usepackage{graphicx}
\DeclareGraphicsExtensions{.eps}
\usepackage{subcaption}

\usepackage{multicol}
\usepackage{multirow}

\usepackage{calc}
\usepackage{tikz}
\usetikzlibrary{shapes, decorations.markings, calc}

\newtheoremstyle{mythmstyle}
  {\topsep}
  {\topsep}
  {}
  {}
  {\bfseries}
  {.}
  {0.5em}
  {}

\theoremstyle{mythmstyle}
\newtheorem{Theorem}{Theorem}

\newtheorem{Corollary}{Corollary}[Theorem]
\newtheorem{Lemma}[Theorem]{Lemma}

\theoremstyle{remark}

\def \ba{\begin{array}}
\def \ea{\end{array}}

\def \bea{\begin{eqnarray}}
\def \eea{\end{eqnarray}}

\def \be{\begin{equation}}
\def \ee{\end{equation}}

\def \BEA{\begin{eqnarray*}}
\def \EEA{\end{eqnarray*}}

\def \BE{\begin{equation*}}
\def \EE{\end{equation*}}

\def \+{\dagger}
\def \mb{\mathbf}
\def \bb{\mathbb}

\def \sf{\mathsf}
\def \mc{\mathcal}

\def \diag{\text{diag}}

\def \sym{\text{sym}}
\def \t{\mathtt{T}}

\def \S{\mathtt{S}}
\def \I{\mathtt{I}}
\def \D{\mathtt{D}}
\def \A{\mathtt{A}}
\def \R{\mathtt{R}}
\def \T{\mathtt{T}}
\def \H{\mathtt{H}}
\def \E{\mathtt{E}}
\def \V{\mathtt{V}}

\def \colsep{\arraycolsep}

\def \ol{\overline}
\def \ul{\underline}

\def\BibTeX{{\rm B\kern-.05em{\sc i\kern-.025em b}\kern-.08em
    T\kern-.1667em\lower.7ex\hbox{E}\kern-.125emX}}

\begin{document}

\title{Observer Design for the State Estimation of Epidemic Processes}

\author{Muhammad Umar B. Niazi
\and
Karl Henrik Johansson
\thanks{The authors are with the Division of Decision and Control Systems and Digital Futures, EECS, KTH Royal Institute of Technology, SE-100 44 Stockholm, Sweden (Email: {\tt mubniazi@kth.se},
{\tt kallej@kth.se}).}
}

\maketitle

\begin{abstract}
Although an appropriate choice of measured state variables may ensure observability, designing state observers for the state estimation of epidemic models remains a challenging task. Epidemic spread is a nonlinear process, often modeled as the law of mass action, which is of a quadratic form; thus, on a compact domain, its Lipschitz constant turns out to be local and relatively large, which renders the Lipschitz-based design criteria of existing observer architectures infeasible. In this paper, a novel observer architecture is proposed for the state estimation of a class of nonlinear systems that encompasses the deterministic epidemic models. The proposed observer offers extra leverage to reduce the influence of nonlinearity in the estimation error dynamics, which is not possible in other Luenberger-like observers.  Algebraic Riccati inequalities are derived as sufficient conditions for the asymptotic convergence of the estimation error to zero under local Lipschitz and generalized Lipschitz assumptions. Equivalent linear matrix inequality formulations of the algebraic Riccati inequalities are also provided. The efficacy of the proposed observer design is illustrated by its application on the celebrated SIDARTHE-V epidemic model.
\end{abstract}


\section{Introduction}
Epidemic processes comprise spreading phenomena in biological, social, geographic, and cyber domains. Therefore, in addition to epidemiology \cite{kiss2017}, modeling and analysis of epidemic processes have found applications in wildfires \cite{somers2020}, cyber-physical systems \cite{tuyishimire2020}, computer networks \cite{mishra2014}, and wireless communications \cite{ojha2019}.
State estimation is crucial for efficient monitoring and control of these processes. Because of the physical and economic limitations of available sensing resources, it is generally not possible to measure all the state variables of an epidemic process. This limitation is well-understood in real-world epidemics, where it is impossible to obtain data on epidemic variables like the number of susceptible, undetected infected, and recovered cases.

Observers are dynamical systems that utilize the available online measurements to estimate the state of a linear/nonlinear process. For this, they require a suitable model of the process. Therefore, before state estimation, it is necessary to identify the model parameters. The notions of identifiability and observability, which have been extensively studied for epidemic models \cite{massonis2021, sauer2021}, ensure whether the parameters can be identified and the state variables can be reconstructed,  respectively. Several parameter estimation techniques for epidemic models are proposed in \cite{hadeler2011, chowell2017, magal2018}. However, despite being equally important, state estimation of epidemic models didn't receive enough attention in the research literature. The existing state estimation techniques for epidemic models either consider a linearized epidemic model \cite{alonso2012} or design an extended Kalman filter \cite{rajaei2021, gomez2021, azimi2022}, which is also a sophisticated estimation technique based on linearization  \cite{bonnabel2014}. The problem with these techniques is that they only guarantee local convergence of the state estimation algorithm. The adaptive observer proposed in \cite{bliman2018} ensures a global convergence result; however, it assumes that the nonlinearity in the model is a measured output, which is a very restrictive assumption.

In this paper, we propose a novel observer for a class of nonlinear systems that encompasses both compartmental and networked epidemic models. The form of the proposed observer combines ideas of feedforward output injection and innovation term in the nonlinearity, which were originally proposed by Luenberger \cite{luenberger1964, luenberger1971} for linear systems and by Arcak and Kokotovi\'{c} \cite{arcak2001} for nonlinear systems, respectively. Using the Lyapunov stability analysis and local Lipschitz property of epidemic processes, we provide a sufficient condition for the asymptotic convergence of estimation error to zero in terms of an algebraic Riccati inequality (ARI), which can be equivalently formulated as a linear matrix inequality (LMI) feasibility condition. Under the assumption of generalized Lipschitz property, we provide a less restrictive ARI as a sufficient condition of the error convergence. The proposed observer is easy to design and provides extra leverage to minimize the influence of nonlinearity in the estimation error dynamics. Such a leverage is not possible in the Luenberger-like observer \cite{thau1973, raghavan1994, rajamani1998} and cannot be fully exploited in the Arcak-Kokotovi\'{c} observer \cite{arcak2001}.

We use standard notations throughout the paper; nonetheless, we clarify wherever an ambiguity can arise. The dependence of signals on time variable $t$ is sometimes omitted for brevity, and it should be considered as implicit. The rest of the paper is organized as follows: In Section~\ref{sec_motivation}, we describe motivation, problems, and challenges in epidemic monitoring and control. After stating the problem in Section~\ref{sec_problem}, we briefly review the existing observer design techniques for nonlinear systems and highlight difficulties in using them for epidemic processes in Section~\ref{sec_literature_observer}. The proposed observer design is presented in Section~\ref{sec_observer}. Finally, Section~\ref{sec_simulation} and \ref{sec_conclusion} provide simulation results and concluding remarks, respectively.

\section{Motivation and Challenges} \label{sec_motivation}

In this section, we define a class of nonlinear systems encompassing the deterministic epidemic models. We emphasize the need to perform identifiability and observability analysis of epidemic models, and choosing the appropriate measured output variables. After briefly reviewing the parameter and state estimation problems, and demarcating the scope of this paper, we list some challenges that should be addressed in future for effective epidemic monitoring and control.

\subsection{Models of epidemic processes}
Epidemic processes are usually described by either population/compartmental or metapopulation/networked models. Compartmental models consider a lumped population divided into multiple epidemic compartments such as susceptible, exposed, infected, recovered, etc. These models assume a homogeneous population structure without an underlying network; therefore, each person is assumed to be equally likely to get infected when the epidemic starts. The basic compartmental models, among others \cite{kiss2017}, are SIS (Susceptible, Infected, Susceptible), SIR (Susceptible, Infected, Removed), and SEIR (Susceptible, Exposed, Infected, Removed). Their variations include SIQR (Susceptible, Infected, Quarantined, Removed) \cite{hamelin2021}, SIDUR (Susceptible, undetected Infected, Detected infected, Unidentified recovered, identified Removed) \cite{niazi2021arc}, SIDARTHE (Susceptible, Infected, Diagnosed, Ailing, Recognized, Threatened, Healed, Extinct) \cite{giordano2020}, and many others (see \cite{massonis2021}).

Networked models are more sophisticated than the compartmental models as they also take into account the geographic and/or demographic heterogeneity of an epidemic process. The underlying network governs the infection transmission from one node to another, \cite{nowzari2016, mei2017, pare2020b}, where the edges between nodes represent interactions between different population classes separated from each other either demographically (e.g., age groups) or geographically. 

All the epidemic models mentioned above belong to a specific class of nonlinear systems described by
\be \label{eq:sys}
\dot{x} = Ax + Gf(Hx), \quad y = Cx
\ee
where $x(t)\in\mc{X}\subset\bb{R}^{n_x}$ is the state and $y(t)\in\bb{R}^{n_y}$ is the measured output. The matrices $A\in\bb{R}^{n_x\times n_x}$, $G\in\bb{R}^{n_x\times n_f}$, $H\in\bb{R}^{k\times n_x}$, and $C\in\bb{R}^{n_y\times n_x}$. The nonlinearity $f:\mc{X}\rightarrow\bb{R}^{n_f}$ is of polynomial type; therefore, it is continuously differentiable and Lipschitz continuous on a bounded domain $\mc{X}$. In other words, for every $x,\hat{x}\in\mc{X}$,
\be \label{eq:lipschitz}
\|f(Hx)-f(H\hat{x})\| \leq \ell \| H(x-\hat{x})\|
\ee
where $\|\cdot\|$ denotes the Euclidean norm and
\be \label{eq:lipschitz_compute}
\ell = \sup_{x\in\mc{X}} \sigma_{\max}\left(\frac{\partial f}{\partial x}(Hx)\right)
\ee
is the Lipschitz constant (see \cite{giaquinta2009}) with $\sigma_{\max}$ the maximum singular value. The problem \eqref{eq:lipschitz_compute} can be solved numerically using a suitable solver for nonlinear constrained optimization problems (e.g., \texttt{fmincon} in MATLAB).

Both SIDARTHE and SIDARTHE-V models, for instance, developed in \cite{giordano2020} and \cite{giordano2021}, respectively, can be written as \eqref{eq:sys} with $f(Hx)=[\colsep=3pt\ba{cccc} \S\I & \S\D & \S\A & \S\R \ea]^\t$ (see Section~\ref{sec_simulation}). Networked SIS epidemic model \cite{pare2020b} in the form of \eqref{eq:sys} has $f(Hx) = -\diag(x)BWx$, $G=H=I_n$, and $A=(BW-D)$, where $B=\diag(\beta_1,\dots,\beta_n)$ is the matrix of infection susceptiblities of nodes, $D=\diag(\delta_1,\dots,\delta_n)$ is the matrix of healing rates of nodes, and $W$ is the weighted adjacency matrix. 

When selecting a suitable model describing an epidemic process, particular consideration should be given to its identifiability and observability properties with respect to the available measurements. The measured output in compartmental models comprise the variables corresponding to the measurable compartments---for example, active number of diagnosed cases and total number of deaths---and/or the variables corresponding to the flows from one compartment to another---for example, daily number of diagnosed and hospitalized cases. Similarly, the measured output in networked models can be obtained, for instance, by recording the active number of diagnosed cases in each node. In the engineering applications \cite{tuyishimire2020, mishra2014, ojha2019}, one can also assume that the full states of some gateway nodes \cite{bullo2022} in the network are measured directly by the dedicated sensors.

\subsection{Identifiability and parameter estimation}

Let $\theta\in\Theta\subset\bb{R}^{n_\theta}$ be the vector of parameters of \eqref{eq:sys}, where $A:=A(\theta)$, $G:=G(\theta)$, and $C:=C(\theta)$. Then, the notion of (global) identifiability ensures that, for any generic $\theta,\hat{\theta}\in\Theta$,
\[
y(t;\theta) \equiv y(t;\hat{\theta}) \implies \theta=\hat{\theta}.
\]
If this implication holds locally within an open neighborhood of every $\theta\in\Theta$, then \eqref{eq:sys} is said to be locally identifiable \cite{audoly2001, saccomani2003}. In other words, for \eqref{eq:sys} to be identifiable, there must exist a time $t>0$ such that the map $\theta\mapsto y(t;\theta)$ is injective, at least locally, for any $\theta\in\Theta$.

There are several softwares like DAISY \cite{bellu2007} and GenSSI \cite{chis2011} that employ tools from differential geometry and differential algebra for checking the identifiability of a nonlinear model from its output $y(t)$ and its derivatives. Identifiability is necessary for the identification of $\theta$ by using the knowledge of output trajectory of the system \eqref{eq:sys}. See \cite{schoukens2019} for an extensive survey on the identification of nonlinear systems.

Identification of epidemic models is crucial for understanding and forecasting the epidemic evolution. It is extensively studied in the literature (see, e.g., \cite{hadeler2011, chowell2017, magal2018}). Given that the system \eqref{eq:sys} is identifiable and that output data is available, the parameter estimation problem can be formulated as a maximum likelihood estimation problem
\be \label{prob:mle_parameters}
\hat{\theta} = \arg\min_{\theta\in\Theta} \sum_{k=0}^{T-1} \|\bar{y}(k) - y(k;\theta) \|^2
\ee
where $\hat{\theta}$ is the estimated parameter vector, $\bar{y}(k)$ is the measurement data, $y(k;\theta)$ is the sampled output from the model \eqref{eq:sys}, and $T$ is the total time duration. If the parameters are time varying, it can be assumed that they are piece-wise constant, as in \cite{giordano2021}. Then, a moving horizon estimation technique \cite{kuhl2011} can be employed to find the optimal parameters. In the current paper, we focus on the state estimation problem of epidemic models by assuming that the parameters $\theta$ have already been estimated.

\subsection{Observability, detectability, and state estimation} \label{subsec_obsdet}

To assess and predict the epidemic evolution and the attained population-level immunity, it is crucial to estimate the unmeasurable states of an epidemic process---for example, number of susceptible, undetected infected, and immune cases. Such information is also necessary to devise optimal epidemic mitigation policies \cite{alamo2021}. 
 
Given the parameters $\theta$, the notion of observability guarantees whether or not the state trajectory $x(t;x_0)$ of \eqref{eq:sys} with $x_0=x(0)$ can be uniquely determined from the output trajectory $y(t; x_0)$. That is, \eqref{eq:sys} is (globally) observable if, for any $x_0,\hat{x}_0\in\mc{X}$,
\be \label{eq:observable}
y(t; x_0) \equiv y(t; \hat{x}_0) \implies x_0 = \hat{x}_0.
\ee
If the implication \eqref{eq:observable} holds locally within an open neighborhood of every $x_0\in\mc{X}$, then \eqref{eq:sys} is said to be locally observable. See \cite{hermann1977} or \cite[Chapter 1]{boutat2021} for all the characterizations of observability and \cite{massonis2021} for the observability analysis of several compartmental epidemic models.


A more flexible notion than observability is detectability, which ensures the following implication
\[
y(t; x_0) \equiv y(t; \hat{x}_0) \implies \lim_{t\rightarrow\infty} \|x(t; x_0)-x(t; \hat{x}_0)\| = 0.
\]
Detectability of \eqref{eq:sys} is a necessary condition for the existence of an (asymptotic) observer, \cite{bernard2022}, which is a dynamical system
\be \label{eq:observer}
\dot{z} = \phi(z,y), \quad \hat{x} = \psi(z,y)
\ee 
with $\phi:\bb{R}^{n_z}\times\bb{R}^{n_y}\rightarrow\bb{R}^{n_z}$ and $\psi:\bb{R}^{n_z}\times\bb{R}^{n_y}\rightarrow\bb{R}^{n_x}$ designed to ensure
$
\lim_{t\rightarrow\infty} \|x(t)-\hat{x}(t)\| = 0.
$
On the other hand, observability of \eqref{eq:sys} is a necessary condition for the existence of a tunable observer, \cite{besancon2007}, where an observer \eqref{eq:observer} is said to be tunable if $\phi,\psi$ can be designed such that, for any $\epsilon>0$, there exists $T>0$ ensuring
$
\|x(t)-\hat{x}(t)\|\leq \epsilon$, $\forall t\geq T.
$
For now, we limit our investigation to asymptotic observers.

\subsection{Challenges in epidemic monitoring and control}

It has been highlighted in \cite{alamo2021} that uncertainties and time delays in data acquisition and reporting pose a big challenge in efficient epidemic monitoring and control. Moreover, epidemic parameters vary with the external conditions such as epidemic mitigation policies, social behavior, and disease mutations, which are highly uncertain and create difficulties in the analysis of epidemic processes. The treatment of such uncertainties and time delays is beyond the scope of this paper. However, we acknowledge that it is a very important research problem and will surely be addressed in our future work. The current paper can be considered as a preliminary step towards the observer design of epidemic processes.

Here, we focus on the technical challenges posed by the models of epidemic processes in the absence of uncertainties and noise. As we saw earlier that epidemic spread is a nonlinear process with quadratic type nonlinearity. Generally speaking, such nonlinearities are only locally Lipschitz and the Lipschitz constant may turn out be relatively large. Moreover, the matrix $A$ representing the linear part of an epidemic process is usually sparse, which makes the LMI-based observer design task difficult. 

\section{Problem Statement} \label{sec_problem}

We assume that the system \eqref{eq:sys} is detectable in a sense defined in Section~\ref{subsec_obsdet}, which is necessary for the existence of $\phi,\psi$ such that an asymptotic observer \eqref{eq:observer} converges, i.e., 
\[
\lim_{t\rightarrow\infty} \|x(t;x_0) - \hat{x}(t;\psi(z_0,Cx_0))\|=0.
\]
Note that if \eqref{eq:sys} is observable, differentially observable, or differentially detectable, then it is also detectable. Thus, detectability is the least restrictive assumption to ensure the well-posedness of the observation problem.

Secondly, we assume that the pair $(A,C)$ is a detectable pair. That is, if there exist $v\in\bb{R}^{n_x}$ and $\lambda\in\bb{C}$ such that $Av = \lambda v$ and $Cv = 0$, then the real part $\text{Re}(\lambda)<0$.

Under the above assumptions, we aim to design an observer of the form \eqref{eq:observer} that asymptotically estimates the state of \eqref{eq:sys}, i.e.,
$
\lim_{t\rightarrow\infty} \|x(t)-\hat{x}(t)\| = 0,
$
given the model $A,G,H,C$ and the output measurements $y(t)$.

\section{Existing Observer Design Techniques} \label{sec_literature_observer}

We briefly review observer design techniques for a class of nonlinear systems \eqref{eq:sys} and discuss their scope and limitations.

\subsection{Luenberger-like observer}
Observer design of nonlinear systems belonging to the class \eqref{eq:sys} is a classical problem in control theory. In this regard, \cite{thau1973} was the first to study the problem and proposed a Luenberger-like observer 
\be \label{eq:luenberger_observer}
\dot{\hat{x}} = A\hat{x} + G f(H\hat{x}) + L (y-C\hat{x})
\ee
where $L\in\bb{R}^{n_x\times n_y}$ is the gain matrix to be designed to ensure that the estimation error $e(t)=x(t)-\hat{x}(t)$ satisfying
\[
\dot{e} = (A-LC)e + G [f(Hx)-f(H\hat{x})]
\]
asymptotically converges to zero.
\subsubsection{Design using Lipschitz property}
The design criteria proposed in \cite{thau1973} was very conservative and worked for a very small subset of applications. Subsequently, multiple approaches \cite{raghavan1994, rajamani1998, zemouche2013} attempted to reduce the conservativeness. However, all these methods, in one way or other, require that the Lipschitz constant $\ell$ is sufficiently small. As the nonlinearity in epidemic models is quadratic in nature, the Lipschitz constant \eqref{eq:lipschitz_compute} may turn out to be large. Thus, Lipschitz-based design criteria of Luenberger-like observers turn out to be infeasible for epidemic models.

\subsubsection{Design using bounded Jacobian property}

As the nonlinearity $f$ is differentiable on the state space $\mc{X}$, it has a bounded Jacobian. A design criterion using bounded Jacobian property is proposed in \cite{phanomchoeng2011} for Luenberger-like observer \eqref{eq:luenberger_observer}. This criterion is based on the modified mean value theorem (MVT) for a vector-valued function $f:\mc{X}\rightarrow \bb{R}^{n_f}$, which postulates that, for any $x,\hat{x}\in\mc{X}$,
\be \label{eq:wrongMVT2}
f(Hx) - f(H\hat{x}) = (\ol{\Delta}\circ\ol{D} + \ul{\Delta}\circ\ul{D})(Hx-H\hat{x})
\ee
where $\circ$ denotes the Hadamard product and matrices $\ol{\Delta}=[\ol{\delta}_{ij}]$, $\ol{D}=[\ol{d}_{ij}]$, $\ul{\Delta}=[\ul{\delta}_{ij}]$, $\ul{D}=[\ul{d}_{ij}]$ are defined such that $0\leq \ol{\delta}_{ij},\ul{\delta}_{ij}\leq 1$ with $\ol{\delta}_{ij}+\ul{\delta}_{ij}=1$, and $\ol{d}_{ij}\geq \max(\partial f_i / \partial x_j)$ and $\ul{d}_{ij}\leq \min(\partial f_i / \partial x_j)$. The proof of this modified MVT in \cite{phanomchoeng2011} assumes an extended version of MVT saying that, for some $c\in\mc{X}$, it holds
\be \label{eq:wrongMVT}
f(Hx) - f(H\hat{x}) = \mb{D}f(Hc) (Hx-H\hat{x})
\ee
where $\mb{D}$ is the Jacobian operator. However, \eqref{eq:wrongMVT} is well-known to be an incorrect extension of MVT to vector-valued functions. A counter example is $f(Hx)=[\ba{cc} \cos(x_1) & \sin(x_1) \ea]^\t$ for $x=[\ba{cc} x_1 & x_2 \ea]^\t$ and $\hat{x}=[\ba{cc} \hat{x}_1 & \hat{x}_2\ea]^\t$ with $x_i,\hat{x}_i\in\bb{R}$. Now, for $x_1=2\pi$, $\hat{x}_1=0$, and any $x_2,\hat{x}_2\in\bb{R}$, we have the left-hand side of \eqref{eq:wrongMVT} equal to zero; however, the right-hand side is non-zero with $Hx-H\hat{x}\neq 0$ for $x-\hat{x}\notin\ker(H)$ and
\[
\mb{D}f(Hc) = \left[\ba{cc} -\sin(c_1) & 0 \\ 0 & \cos(c_1) \ea\right] \neq 0, \, \forall\,c\in\bb{R}^2.
\]
Therefore, the results \eqref{eq:wrongMVT2} and \eqref{eq:wrongMVT} claimed in \cite{phanomchoeng2011} do not hold in general. See \cite[Section 1.3.2 (b)]{giaquinta2009} for a correct extension of MVT, which is in an integral form.

\subsubsection{Design using one-sided Lipschitz and quadratically inner boundedness properties}

A function $f:\mc{X}\rightarrow\bb{R}^{n_f}$ is said to be one-sided Lipschitz (OSL) if, for every $x,\hat{x}\in\mc{X}$, there exists $\varphi\in\bb{R}$ such that
\be \label{eq:osl}
\langle G(f(Hx)-f(H\hat{x})), x-\hat{x} \rangle \leq \varphi \| x-\hat{x}\|^2 .
\ee
It is said to be quadratically inner bounded (QIB) if, for every $x,\hat{x}\in\mc{X}$, there exist $\varrho_1,\varrho_2\in\bb{R}$ such that
\be \label{eq:qib}
\ba{l}
\|G(f(Hx)-f(H\hat{x}))\|^2 \leq \varrho_1 \|x-\hat{x}\|^2  \\ \qquad\qquad + \varrho_2 \langle G(f(Hx)-f(H\hat{x})), x-\hat{x} \rangle.
\ea\ee
See \cite{nugroho2021} for methods to compute the constants $\varphi,\varrho_1,\varrho_2$.

The design criteria of a Luenberger-like observer \eqref{eq:luenberger_observer} using the OSL and QIB properties was first provided in \cite{abbaszadeh2010}. However, the design is still conservative as the one-sided Lipschitz property cannot be used directly. To elucidate, consider a Lyapunov function $V=e^\t P e$, where $P=P^\t>0$, then
\[
\dot{V} = e^\t [(A-LC)^\t P + P(A-LC)] e + 2e^\t P G\tilde{f}
\]
where $\tilde{f}:=f(Hx)-f(H\hat{x})$.
Thus, OSL property must hold for $e^\t P G\tilde{f}$, instead of $e^\t G\tilde{f}$ as in \eqref{eq:osl}, for some $P>0$. Assuming $P=I$ significantly, and unnecessarily, increases the conservativeness. The technique proposed by \cite{abbaszadeh2010} tries to reduce the difference between the maximum and minimum eigenvalues of $P$. Another technique proposed by \cite{zhao2010} constrains $P$ to be block diagonal to satisfy OSL property \eqref{eq:osl} for $Pe$ instead of $e=x-\hat{x}$. However, these assumptions hold for very peculiar nonlinearities and turn out to be quite conservative in general. An approach based on Algebraic Riccati Equation presented in \cite{zhang2012} results in a Linear Matrix Inequality, which uses the OSL property in the QIB \eqref{eq:qib} and obtain a Lipschitz-type inequality. This method turns out to be as restrictive as those based on the Lipschitz property.

\subsection{Arcak-Kokotovi\'{c} observer}

An interesting approach proposed by Arcak and Kokotovi\'{c} \cite{arcak2001} considers the following observer:
\be \label{eq:arcak_kokotovic}
\dot{\hat{x}} = A\hat{x} + Gf(H\hat{x}+K(y-C\hat{x})) + L(y-C\hat{x})
\ee
where $K\in\bb{R}^{k\times n_y}$ and $L\in\bb{R}^{n_x\times n_y}$ are two  design matrices. The estimation error $e(t)=x(t)-\hat{x}(t)$ satisfies
\[
\dot{e} = (A-LC)e + G \tilde{f}(x,\hat{x},y)
\]
where  
$
\tilde{f}(x,\hat{x},y) := f(Hx)-f(H\hat{x}+K(y-C\hat{x})).
$
The Lipschitz property \eqref{eq:lipschitz} results in
\be \label{eq:lipschitz_2}
\|\tilde{f}\| \leq \ell \|(H-KC)e\|.
\ee
Thus, we have an extra leverage to minimize the right-hand side of \eqref{eq:lipschitz_2} by choosing appropriate $K$. The design proposed in \cite{arcak2001} is catered to a special class of systems where $f(Hx)$ is nondecreasing, which is not the case with epidemic models. This limitation is removed in \cite{zemouche2017}, which provides an equivalent characterization of the Lipschitz property \eqref{eq:lipschitz_2}. However, this characterization is time-varying and is not easy to satisfy in general. Finally, \cite{rajamani2020} designs the observer \eqref{eq:arcak_kokotovic} with switched gains; however, designing the switching signal remains an open problem.

\subsection{Observer design using nonlinear state transformation}

Assume there exists an injective map $T:\bb{R}^{n_x}\rightarrow\bb{R}^{n_z}$, where $n_z \geq n_x$, such that the coordinate transformation $z=T(x)$ of \eqref{eq:sys} yields another system in $z(t)\in\bb{R}^{n_z}$, which has a specific structure (triangular, normal, etc.) or it is linear upto output injection. Then, several observer design methods like high-gain, backstepping, finite-time, and Kazantis-Kravaris-Luenberger become convenient. See \cite{besancon2007, boutat2021, bernard2022} for more details on these observer design techniques. However, the problem with this approach is that, once the estimate $\hat{z}(t)$ is obtained, one needs to find the inverse of the transformation $T$ to obtain the estimate $\hat{x}(t)$ in the original coordinates. Unless $T$ is a diffeomorphism, which means $n_z=n_x$, obtaining a left inverse of $T$ is very challenging \cite{andrieu2021}. The papers \cite{ramos2020, peralez2021} approximate $T$ and its left inverse by a neural network. However, the neural network results in overfitting and doesn't generalize well when the real trajectory is significantly different from the training trajectories.

\section{Proposed Observer Design} \label{sec_observer}

We propose an observer of the form
\be \label{eq:our_observer}
\ba{ccl}
\dot{z} &=& Mz + (ML+J)y + N Gf(v) \\
v &=& H\hat{x} + K (y-C\hat{x}) \\
\hat{x} &=& z + Ly
\ea
\ee
where
\be \label{eq:our_matrices}
M = A-LCA-JC, \quad N = I-LC
\ee
and $J,L\in\bb{R}^{n_x\times n_y}$ and $K\in\bb{R}^{k\times n_y}$ are gain matrices to be designed.
The estimation error $e(t)=x(t)-\hat{x}(t)=Nx(t) - z(t)$ satisfies
\be \label{eq:our_error}
\dot{e} = M e + N G [f(Hx) - f(v)]
\ee
where $M,N$ are given in \eqref{eq:our_matrices}.

The form of the observer \eqref{eq:our_observer} enables one to design $L$ and $K$ so that the influence of $NG [f(Hx) - f(v)]$ in the error dynamics \eqref{eq:our_error} is filtered out. This is made possible by appropriately choosing the feedforward output injection $Ly$ in \eqref{eq:our_observer}, which is inspired by the Luenberger observer \cite{luenberger1964, luenberger1971} for linear systems and Daruoach observer \cite{darouach1994} for linear systems with unknown inputs. Also, the internal innovation term $K(y-C\hat{x})$ in the function $f$, which is inspired by the Arcak-Kokotovic observer \eqref{eq:arcak_kokotovic} for monotonic nonlinear systems (see \cite{arcak2001}), enables one to choose $K$ such that the difference $f(Hx)-f(v)$ in the error dynamics is minimized. These leverages allow to reduce the conservativeness of the Lipschitz-based design criteria of our observer, which is not possible in the contemporary observers presented in Section~\ref{sec_literature_observer}. Finally, note that $J$ can always be chosen such that $M=A-LCA-JC$ is Hurwitz given that $(A,C)$ is a detectable pair.

\subsection{Design criterion using Lipschitz property}

We present a design criterion for the gain matrices $J$, $L$, and $K$ based on the Lipschitz property \eqref{eq:lipschitz_2} to guarantee the asymptotic stability of the (estimation) error dynamics \eqref{eq:our_error}.

\begin{Theorem} \label{thm:ARI1}
Subject to the Lipschitz property \eqref{eq:lipschitz_2}, if there exist a positive definite matrix $P\in\bb{R}^{n_x\times n_x}$ and design matrices $J,L\in\bb{R}^{n_x\times n_y}$ and $K\in\bb{R}^{k\times n_y}$ such that the following algebraic Riccati inequality holds
\be \label{eq:ari}
\ba{l}
(A-LCA-JC)^\t P + P (A-LCA-JC) \\
\qquad +\, P (I-LC) GG^\t (I-LC)^\t P \\
\qquad\qquad +\, \ell^2 (H-KC)^\t (H-KC) < 0
\ea
\ee
then the error dynamics \eqref{eq:our_error} is asymptotically stable. 
\end{Theorem}
\begin{IEEEproof}
For the asymptotic stability of the error dynamics \eqref{eq:our_error}, we consider a radially unbounded, positive definite function $V:\bb{R}^{n_x}\setminus\{0\}\rightarrow\bb{R}_{>0}$ with $V(0)=0$, and show that it is a Lyapunov function, i.e., $\dot{V}<0$, if \eqref{eq:ari} is satisfied. Let $V(e(t))=e^\t(t) P e(t)$ be such a Lyapunov function with $P>0$. Then, by differentiating $V$ with respect to time $t$ along the trajectories of the error dynamics \eqref{eq:our_error}, we obtain
\be \label{eq:V_dot}
\colsep=1pt\ba{ccl}
\dot{V} &=& \dot{e}^\t P e + e^\t P \dot{e} \\
&=&  e^\t (M^\t P + P M) e + 2e^\t P N G \tilde{f}
\ea\ee
where $\tilde{f} := f(Hx) - f(H\hat{x} + K (y-C\hat{x}))$. Notice that
\[
\colsep=2pt\ba{ccl}
2e^\t P N G \tilde{f} &\leq& |2 e^\t P N G\tilde{f}| \\
&\leq & 2 \|G^\t N^\t P e\| \|\tilde{f}\| \\ &\leq & \| G^\t N^\t P e \|^2 + \|\tilde{f}\|^2
\ea\]
where the last two steps are due to Cauchy-Schwarz and Young's inequalities, respectively.
Then, from \eqref{eq:lipschitz_2}, we have
\[
\|\tilde{f}\|^2 \leq \ell^2 \| (H-KC) e\|^2
\]
implying
\[
2e^\t P N G \tilde{f} \leq e^\t P N GG^\t N^\t P e + \ell^2 e^\t (H-KC)^\t (H-KC) e.
\]
Using it in \eqref{eq:V_dot} and using the values of $M$ and $N$ from \eqref{eq:our_matrices}, we obtain that if \eqref{eq:ari} holds then $\dot{V}<0$ for every $e\neq 0$.
\end{IEEEproof}

The ARI in \eqref{eq:ari} can be equivalently formulated as an LMI condition, whose feasibility can be checked using a suitable semidefinite programming software (e.g., YALMIP).

\begin{Corollary} \label{cor_LMI1}
The ARI \eqref{eq:ari} holds if and only if there exist positive definite matrices $P,Q\in\bb{R}^{n_x\times n_x}$, $K\in\bb{R}^{k\times n_y}$, and $R,S\in\bb{R}^{n_x\times n_y}$ such that the following LMIs are feasible
\begin{subequations} \label{eq:lmi}
\begin{align}
\left[\ba{cc}
\sym(PA - RCA - SC) + Q & (P-RC)G \\
G^\t(P-RC)^\t & -I_{n_f}
\ea\right] &< 0  \label{eq:lmi1} \\
\left[\ba{cc}
-Q & (H-KC)^\t \\
H-KC & -\frac{1}{\ell^2} I_k
\ea\right] &\leq 0 \label{eq:lmi2}
\end{align}
\end{subequations}
where $\sym(X)=X+X^\t$, and $J=P^{-1} S$ and $L=P^{-1} R$.
\end{Corollary}
\begin{IEEEproof}
Assume that \eqref{eq:lmi} holds. Then, by Schur complement lemma \cite[Chapter 1]{boyd1994}, \eqref{eq:lmi2} is equivalent to
\[
\ell^2 (H-KC)^\t (H-KC) \leq Q.
\]
Thus, if
\[\ba{r}
(A-LCA-JC)^\t P + P(A-LCA-JC) \qquad \\ + P(I-LC)GG^\t(I-LC)^\t P + Q <0
\ea\]
then \eqref{eq:ari} holds. Substituting $R=PL$ and $S=PJ$, we can write
\be \label{eq:ari_3}
\sym(PA-RCA-SC) + (P-RC) GG^\t (P-RC)^\t + Q < 0.
\ee
Again, by Schur complement lemma, \eqref{eq:ari_3} is equivalent to \eqref{eq:lmi1}. Thus, if \eqref{eq:lmi} holds, then \eqref{eq:ari} is satisfied. 

For proving the necessity, assume \eqref{eq:lmi1} holds but not \eqref{eq:lmi2}, then $\ell^2(H-KC)^\t (H-KC)> Q$. Thus, even if \eqref{eq:lmi1} holds, \eqref{eq:ari} cannot be guaranteed. Second, assume \eqref{eq:lmi2} holds but not \eqref{eq:lmi1}, then, even for $\ell^2(H-KC)^\t (H-KC)=Q$, \eqref{eq:ari_3} doesn't hold. Thus, \eqref{eq:ari} doesn't hold either.
\end{IEEEproof}

The LMI condition \eqref{eq:lmi} is less restrictive than the design criteria of Luenberger-like and Arcak-Kokotovi\'{c} observers presented in the previous section. However, using the generalized Lipschitz condition, we can further reduce the conservativeness of the above design criterion.

\subsection{Design criterion using generalized Lipschitz property}

Lipschitz property is considered to be conservative, especially when the Lipschitz constant turns out to be very large. The generalized Lipschitz condition, introduced in \cite{ekramian2011}, is considered to be less conservative, which states that, for every $x,\hat{x}\in\mc{X}$, there exist positive definite matrices $V\in\bb{R}^{n_f\times n_f}$ and $W\in\bb{R}^{k\times k}$ such that
\be \label{eq:genLipschitz}
\|f(Hx)-f(H\hat{x})\|_V \leq \|H(x-\hat{x})\|_W
\ee
where the norms $\|a\|_V:=\sqrt{a^\t V a}$ and $\|b\|_W:=\sqrt{b^\t W b}$ for $a\in\bb{R}^{n_f}$ and $b\in\bb{R}^k$.

\begin{Lemma} \label{lemma_genLip}
If the generalized Lipschitz property \eqref{eq:genLipschitz} holds for some positive definite matrices $V$ and $W$, then, for every $\Gamma\in\bb{R}^{n_f\times k}$ and $x,\hat{x}\in\mc{X}$, it holds
\[
2 \xi^\t \Gamma^\t \tilde{f} \leq \|\Gamma \xi\|_{V^{-1}}^2 + \|\xi\|_W^2
\]
where $\xi:=H(x-\hat{x})$ and $\tilde{f}:=f(Hx)-f(H\hat{x})$.
\end{Lemma}
\begin{IEEEproof}
It holds that
\[\ba{rcl}
0 &\leq & (V^{-\frac{1}{2}} \Gamma \xi - V^{\frac{1}{2}} \tilde{f})^\t (V^{-\frac{1}{2}} \Gamma \xi - V^{\frac{1}{2}} \tilde{f}) \\
&=& \|\Gamma\xi\|_{V^{-1}}^2 + \|\tilde{f}\|_V^2 - 2\xi^\t\Gamma^\t \tilde{f} \\
&\leq & \|\Gamma\xi\|_{V^{-1}}^2 + \|\xi\|_W^2 - 2\xi^\t\Gamma^\t \tilde{f}.
\ea\]
\end{IEEEproof}

Similar to Theorem~\ref{thm:ARI1} that employs the Lipschitz property, we present an ARI-based design criterion that employs the generalized Lipschitz property.

\begin{Theorem}
Subject to the generalized Lipschitz condition \eqref{eq:genLipschitz}, if there exist a positive definite matrix $P\in\bb{R}^{n_x\times n_x}$ and design matrices $J,L\in\bb{R}^{n_x\times n_y}$ and $K\in\bb{R}^{k\times n_y}$ such that the following algebraic Riccati inequality holds
\be \label{eq:ari2}
\ba{l}
(A-LCA-JC)^\t P + P (A-LCA-JC) \\
\qquad +\, P (I-LC) G V^{-1} G^\t (I-LC)^\t P \\
\qquad\qquad +\, (H-KC)^\t W (H-KC) < 0
\ea
\ee
then the error dynamics \eqref{eq:our_error} is asymptotically stable. 
\end{Theorem}
\begin{IEEEproof}
By using Lemma~\ref{lemma_genLip} in \eqref{eq:V_dot}, it is straightforward to see that $\dot{V}<0$ if \eqref{eq:ari2} holds.
\end{IEEEproof}

If one can find $V,W$ that satisfy the generalized Lipschitz condition \eqref{eq:genLipschitz}, then it is trivial to see that the ARI \eqref{eq:ari2} is less restrictive than the ARI \eqref{eq:ari}.
As in Corollary~\ref{cor_LMI1}, we can find equivalent LMI condition for the ARI \eqref{eq:ari2} as follows.

\begin{Corollary}
The ARI \eqref{eq:ari2} holds if and only if there exist positive definite matrices $P,Q\in\bb{R}^{n_x\times n_x}$, $K\in\bb{R}^{k\times n_y}$, and $R,S\in\bb{R}^{n_x\times n_y}$ such that the following LMIs are feasible
\begin{subequations} \label{eq:lmi_2}
\begin{align}
\left[\ba{cc}
\sym(PA - RCA - SC) + Q & (P-RC)G \\
G^\t(P-RC)^\t & -V
\ea\right] &< 0  \label{eq:lmi1_2} \\
\left[\ba{cc}
-Q & (H-KC)^\t \\
H-KC & -W^{-1}
\ea\right] &\leq 0 \label{eq:lmi2_2}
\end{align}
\end{subequations}
where $J=P^{-1} S$ and $L=P^{-1} R$.
\end{Corollary}
\begin{IEEEproof}
The proof is similar to that of Corollary~\ref{cor_LMI1}.
\end{IEEEproof}

\section{Simulation Results} \label{sec_simulation}

\begin{figure}[!]
\includegraphics[trim={0 0 0 18}, clip, width=0.5\textwidth]{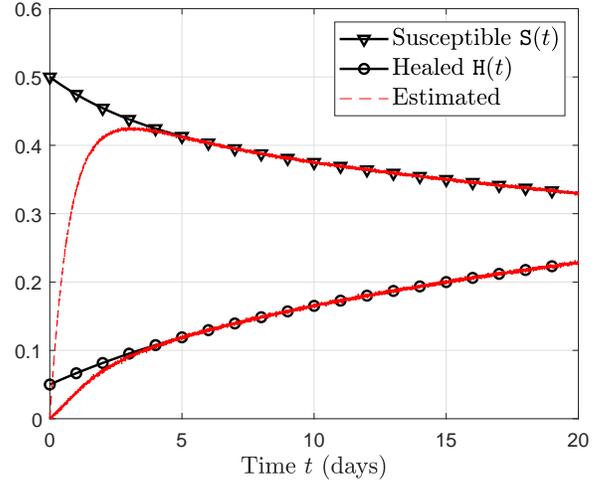}
\caption{Estimation of susceptible and healed cases.}
\label{fig:estimation}
\end{figure}

\begin{figure}[!]
\includegraphics[trim={0 0 0 18}, clip, width=0.5\textwidth]{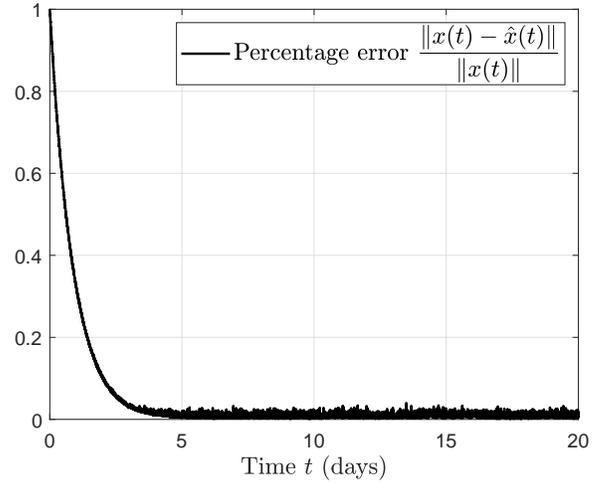}
\caption{Percentage estimation error.}
\label{fig:error}
\end{figure}

For illustrating the proposed observer design, we consider the SIDARTHE-V model \cite{giordano2021} described by
\be \label{eq:sidarthe-v}
\ba{ccl}
\dot{\S} &=& -\S(\alpha \I + \beta \D + \gamma \A + \delta \R) - \varphi \S \\
\dot{\I} &=& \S(\alpha \I + \beta \D + \gamma \A + \delta \R) - (\varepsilon + \zeta + \lambda) \I \\
\dot{\D} &=& \varepsilon \I - (\eta + \rho) \D \\
\dot{\A} &=& \zeta \I - (\theta + \mu + \kappa) \A \\
\dot{\R} &=& \eta \D + \theta \A - (\nu + \xi + \tau_1) \R \\
\dot{\T} &=& \mu \A + \nu \R - (\sigma + \tau_2) \T \\
\dot{\H} &=& \lambda \I + \rho \D + \kappa \A + \xi \R + \sigma \T \\
\dot{\E} &=& \tau_1 \R + \tau_2 \T \\
\dot{\V} &=& \varphi \S
\ea
\ee
where the states correspond to the proportion of Susceptible, Infected (asymptomatic, undetected), Diagnosed (asymptomatic, detected), Ailing (symptomatic, undetected), Recognized (symptomatic, detected), Threatened (acutely symptomatic, detected), Healed (recovered naturally or with treatment), Extinct (died due to the disease). That is, the state space $\mc{X}=[0,1]^{n_x}$ with $n_x=9$. For the definition and interpretation of all the parameters, please refer to \cite{giordano2021}. Below, except $\varphi$, which is chosen arbitrarily, we list the parameters as estimated for the last segment of the time period considered in \cite{giordano2021}:
\[\ba{cccc}
\alpha = 0.3872 & \beta = 0.0053 & \gamma =  0.1485 & \delta = 0.0050 \\
\varepsilon = 0.2988 & \theta = 0.3700 & \zeta = 0.0025 & \eta = 0.0018 \\
\mu = 0.1200 & \nu = 0.0200 & \tau_1 = 0.0050 & \tau_2 = 0.1700 \\
\lambda = 0.1128 & \rho = 0.0320 & \kappa = 0.0200 & \xi = 0.0120 \\
\multicolumn{2}{r}{\sigma = 0.0240} & \multicolumn{2}{l}{\varphi = 0.0500.}
\ea\]
The measured output corresponds to the available data for COVID-19, and is given by
\[
y = \underbrace{\left[\arraycolsep=5pt\ba{ccccccccc}
0 & \varepsilon & 0 & \theta+\mu & 0 & 0 & 0 & 0 & 0 \\
0 & 0 & 1 & 0 & 0 & 0 & 0 & 0 & 0 \\
0 & 0 & 0 & 0 & 1 & 0 & 0 & 0 & 0 \\
0 & 0 & 0 & 0 & 0 & 1 & 0 & 0 & 0 \\
0 & 0 & 0 & 0 & 0 & 0 & 0 & 1 & 0 \\
0 & 0 & 0 & 0 & 0 & 0 & 0 & 0 & 1 \\
1 & 1 & 1 & 1 & 1 & 1 & 1 & 1 & 1 
\ea\right]}_{C} \underbrace{\left[\ba{c} \S \\ [-2.5pt] \I \\ [-2.5pt] \D \\ [-2.5pt] \A \\ [-2.5pt] \R \\ [-2.5pt] \T \\ [-2.5pt] \H \\ [-2.5pt] \E \\ [-2.5pt] \V \ea\right]}_{x}.
\]
The elements of $y$ respectively comprise the normalized number of daily detected cases, active asymptomatic cases, active symptomatic cases, hospitalized cases, deaths, vaccinated cases, and total population. The measured output is assumed to be corrupted by an additive zero-mean white noise vector.

The linear and nonlinear terms in \eqref{eq:sidarthe-v} can be collected in $Ax$ and $Gf(Hx)$ in \eqref{eq:sys}, respectively, where $f(Hx)=[\ba{cccc} \S\I & \S\D & \S\A & \S\R \ea]^\t$ with $H=[\ba{cc} I_5 & 0_{5\times 4} \ea]$ and $A,G$ obtained accordingly from \eqref{eq:sidarthe-v}. The Lipschitz constant \eqref{eq:lipschitz_compute} is obtained to be $\ell=1$. For the observer \eqref{eq:our_observer}, we solve the LMI \eqref{eq:lmi} and obtain a feasible solution $J,L,K$, where
\[\ba{l}
J \hspace*{-2pt}=\hspace*{-3pt} \left[{\scriptsize \colsep=2pt\ba{rrrrrrrrr}
 0.0000 &   0.0000 &   0.0000 &   0.0000 &   0.0000 &   0.0000 &   0.0000 \\
 104.1222 &   0.0000 &   0.0000 &   0.0000 &   0.0000 &   0.0000 &   0.0000 \\
 0.0000 &   0.7771 &   0.0000 &   0.0000 &   0.0000 &   0.0000 &   0.0000 \\
 0.0000 &   0.0000 &   0.0000 &   0.0000 &   0.0000 &   0.0000 &   0.0000 \\
 0.0000 &  -0.0202 &   0.0000 &   0.0000 &   0.0000 &   0.0000 &   0.0000 \\
 0.0000 &   0.0000 &   0.0000 &   0.0000 &   0.0000 &   0.0000 &   0.0000 \\
-40.4200 &   0.0000 &   0.0000 &   0.0000 &  -1.1936 &  -1.1936 &   1.1936 \\
 0.0000 &   0.0000 &   0.0000 &   0.0000 &   1.1508 &   0.0000 &   0.0000 \\
 0.0000 &   0.0000 &   0.0000 &   0.0000 &   0.0000 &   1.1508 &   0.0000
\ea}\right] \\ [3.5em]
L \hspace*{-2pt}=\hspace*{-3pt} \left[{\scriptsize \colsep=0.45pt\ba{rrrrrrrrr}
-3.3467 &  -0.0775 &  -1.4559 &  -2.0955 &  -2.3914 & 119.0938 &   0 \\
 3.3467 & -99.8836 & -95.6678 &-111.7653 &-127.2438 & -46.2617 &   0 \\
 0.0000 &  -9.7224 &   0.0000 &   0.0000 &   0.0000 &   0.0000 &   0 \\
 0.0000 &  -5.1274 & -96.2807 & 302.5176 & 356.3921 &   0.0000 &   0 \\
 0.0000 &  -2.2362 & -29.7965 &   0.0000 &   0.0000 &   0.0000 &   0 \\
 0.0000 &  -0.1708 &  -3.2065 &  -4.9320 &   0.0000 &   0.0000 &   0 \\ 
 0.0000 &  36.7216 &  44.2257 &  19.0293 &  14.8359 &-127.8999 &   0 \\
 0.0000 &   0.0000 &   0.0000 &   0.0000 &   1.0000 &   0.0000 &   0 \\
 0.0000 &   0.0000 &   0.0000 &   0.0000 &   0.0000 &   1.0000 &   0
\ea}\right] \\ [3.5em]
K \hspace*{-4pt}=\hspace*{-3pt} \left[{\scriptsize \colsep=1pt\ba{rrrrrrrrr}
-0.9283 &  -0.3876 &  -0.3876 &  -0.3876 &  -0.3876 &  -0.3876 &   0.3876 \\
 0.4916 &  -0.1103 &  -0.1103 &  -0.1103 &  -0.1103 &  -0.1103 &   0.1103 \\
 0.0000 &   0.8329 &   0.0000 &   0.0000 &   0.0000 &   0.0000 &   0.0000 \\
 1.4001 &   0.0672 &   0.0672 &   0.0672 &   0.0672 &   0.0672 &  -0.0672 \\
 0.0000 &   0.0000 &   0.8329 &   0.0000 &   0.0000 &   0.0000 &   0.0000
\ea}\right].
\ea\]
The initial condition of the observer is chosen to be $z(0)=-Ly(0)$. In Figure~\ref{fig:estimation}, we illustrate the estimation of two unmeasured states, $\S(t)$ and $\H(t)$. The overall estimation performance is illustrated in Figure~\ref{fig:error}, which shows the norm of the estimation error as compared to the norm of the state is about 2\% under noisy measurements. The estimation error is stable under measurement noise and converges to zero asymptotically in the absence of noise.

\section{Conclusion} \label{sec_conclusion}

State estimation of epidemic processes is crucial to assess and predict the infection spread in a population. It is also necessary for devising mitigation policies and optimal allocation of resources. However, it is a challenging problem and the existing observer design techniques turn out to be infeasible for epidemic processes. We proposed a novel observer architecture for a class of nonlinear systems that encompasses both the compartmental and networked epidemic models described by a system of ordinary differential equations. The proposed architecture provides extra leverage against a large Lipschitz constant. This aspect makes the proposed observer less restrictive than other contemporary observers. In addition to the design criterion under Lipschitz property, we proposed another criterion based on generalized Lipschitz condition, which is less restrictive than the Lipschitz condition. Other parametrizations of the system's nonlinearity, like bounded Jacobian or one-sided Lipschitz conditions, can also be considered. The future prospects of the proposed observer architecture include robust design for epidemic processes under model uncertainties and measurement noise.

\section{Acknowledgment}
This work was supported by the Swedish Research Council and the Knut \& Alice Wallenberg Foundation, Sweden.


\end{document}